\documentclass[preprint,12pt]{elsarticle}

    \newcommand{\R}{\mathbb{R}}     
    \def\S{\mathbb{S}}
    
\usepackage{amssymb}
\usepackage{amsmath}
\usepackage{color}

\usepackage[utf8]{inputenc}

\begin{document}

\makeatletter
\def\ps@pprintTitle{%
  \let\@oddhead\@empty
  \let\@evenhead\@empty
  \let\@oddfoot\@empty
  \let\@evenfoot\@oddfoot
}
\makeatother

\begin{frontmatter}
 \newtheorem{thm}{Theorem}[section]
 \newtheorem{lem}[thm]{Lemma} 
  \newtheorem{proposition}[thm]{Proposition}
 \newdefinition{rmk}[thm]{Remark} 
  \newdefinition{cor}[thm]{Corollary} 
 \newproof{pf}{Proof}
  \newdefinition{definition}{Definition}
 
\title{Radial solutions for equations of Weingarten type} 
\author[1]{Antonio Bueno}
\ead{antonio.bueno@cud.upct.es}
\author[2]{Rafael L\'opez} 
\ead{rcamino@ugr.es}
\address[1]{Departamento de Ciencias, Centro Universitario de la Defensa de San Javier. 30729  Santiago de la Ribera, Spain}
\address[2]{Departamento de Geometr\'{\i}a y Topolog\'{\i}a, Universidad de Granada. 18071 Granada, Spain}
\begin{abstract}
In this paper we study the linear Weingarten equation defined by the fully non-linear PDE
 $$a\, \mbox{div}\frac{Du}{\sqrt{1+|Du|^2}}+b\, \frac{\mbox{det}D^2u}{(1+|Du|^2)^2}=\phi\left(\frac{1}{\sqrt{1+|Du|^2}}\right)$$
in a domain $\Omega\subset\R^2$, where $\phi\in C^1([-1,1])$ and $a,b\in\R$. We approach the existence of radial solutions when $\Omega$ is a disk of small radius, giving an affirmative answer when the PDE is of elliptic type. In the hyperbolic case we show that no radial solution exists, while in the parabolic case we find explicitly all the solutions. Finally, in the elliptic case we prove uniqueness and symmetry results concerning the Dirichlet problem of such equation.
\end{abstract}


\end{frontmatter}
 
\section{Introduction}
Consider the existence and uniqueness of classical solutions for the Dirichlet problem
\begin{equation}\label{w1}
\left\{\begin{array}{ll}
a\, \mbox{div} \dfrac{Du}{\sqrt{1+|Du|^2}} +b\, \dfrac{\mbox{det} D^2u}{(1+|Du|^2)^2}=\phi\left(\dfrac{1}{\sqrt{1+|Du|^2}}\right), & \mbox{in $\Omega$}\\
 u=0,& \mbox{on $\partial \Omega$,}
\end{array}\right.
\end{equation}
where $\Omega$ is a bounded smooth domain of $\R^2$, $a,b$ are constants and $\phi\in C^1([-1,1])$. In equation \eqref{w1}, the first term in the left-hand side is a quasilinear operator, while the second one is of  Monge-Amp\`ere type. Some equations of paramount importance have already appeared in the literature for particular choices of the constants $a,b$. For example, if $b=0$, then \eqref{w1} falls in the class of prescribed mean curvature equations where the right-hand side depends of the gradient $Du$. This equation has attracted the attention of many researchers, becoming a fruitful topic of interest. Without aiming to collect all this bibliography, we refer to the reader to \cite{oo} and references therein. 

A solution of \eqref{w1} parametrizes a surface in Euclidean space $\R^3$ whose mean curvature $H$ and Gauss curvature $K$ satisfy the relation 
\begin{equation}\label{w2}
2aH+bK=\phi(\langle N,v\rangle),
\end{equation}
where $N$ is the Gauss map of the surface and $v=(0,0,1)$. In general, a surface that satisfies  a relationship $W(H,K)=0$ between $H$ and $K$ is called a Weingarten surface. The simplest relation $W$ is being linear, that is, $2aH+bK=c$ for constants $a,b,c\in\R$. Regarding this equation, surfaces with constant mean curvature ($b=0$) and with constant Gauss curvature ($a=0$) are particular examples of linear Weingarten surfaces. From now, we will suppose that $a,b\not=0$. 

The generalization  \eqref{w2} is motivated by the theory of the flow by the mean curvature of Huisken, Sinestrari and Ilmanen   \cite{hsi,il} and the flow by the Gauss curvature of Andrews and Urbas (\cite{an,ur}). For example, a  translating soliton  $S$ is  a solution of the mean curvature flow when $S$ evolves purely by translations along some direction of the space. If this direction is, say, $v=(0,0,1)$, then $S+t v$, $t\in \R$, satisfies that for fixed $t$, the normal component of the velocity vector $v$   at each point is equal to the mean curvature at that point. For the initial surface $S$, this implies that $2H=\langle N,v\rangle$. In nonparametric form, $\langle N,v\rangle$ coincides with $1/\sqrt{1+|Du|^2}$, so the surface satisfies \eqref{w1} for $b=0$ and $\phi$ the identity. Similarly, translating solitons by the Gauss curvature are obtained in the same fashion by doing in \eqref{w1} $a=0$ and $\phi$ the identity. Finally, we point out that the first author, together with G\'alvez and Mira, have developed a theory of  complete  surfaces   whose mean curvature is given as a prescribed function of its Gauss map, generalizing some well-known results of the theory of   constant mean curvature surfaces  and translating solitons of the mean curvature flow (\cite{bgm,bgm2}). 

The purpose of this paper is to investigate the radial solutions of \eqref{w1} when $\Omega$ is a round ball $B(0,R)\subset\R^2$ centered at the origin $0$ and of radius $R>0$. It is also desirable that the solutions of \eqref{w1} inherit the symmetries of $\Omega$, so if $\Omega$ is a round ball, a solution of \eqref{w1} must be radial. Our interest is to determine the existence and uniqueness of radial solutions starting orthogonally from the rotation axis. In the case that $u=u(r)$ is such a radial solution, equation \eqref{w1} transforms into the initial value problem 
\begin{equation}\label{rot}
\left\{\begin{array}{ll}a\left(\dfrac{u''}{(1+u'^2)^{3/2}}+\dfrac{u'}{r\sqrt{1+u'^2}}\right)+b\dfrac{u''u'}{r(1+u'^2)^2}=\phi\left(\dfrac{1}{\sqrt{1+u'^2}}\right), & \mbox{in $(0,R)$}\\
 u(0)=0, u'(0)=0,
  \end{array}\right.
\end{equation} 
Let us notice that the equation in  \eqref{rot} is singular at $r=0$, so the solvability is not assured by standard methods. Equivalently, we are asking for the existence of rotational surfaces satisfying the Weingarten relation \eqref{w2} whose generating curve meets orthogonally the rotation axis.   In the field of geometry, there is a great interest in the classification of rotational linear Weingarten surfaces (replacing $\phi$ by a constant $c$) in Euclidean space (\cite{lo1,rs,st}) and also in other ambient spaces (\cite{bss,gm,lo1,lo2,mr}). As a consequence of our investigations, we realized that the existence of solutions of \eqref{rot} depended not only on the constants $a$, $b$ and the function $\phi$, but strongly also on the character of \eqref{rot} as a partial differential equation.  In case that the equation is elliptic at $r=0$, we give a positive answer to the problem. 

\begin{thm} \label{t1}
Suppose that \eqref{rot} is elliptic at $r=0$. Then there is $R>0$ such that   there exists a  solution of \eqref{rot} in $[0,R]$.
\end{thm} 
The elliptic caracter of \eqref{rot} at $r=0$ depends on the sign of $a^2+b\phi(1)$. For the particular case $2aH+bK=c$ and when this relation is elliptic, we provide a proof of the existence of solutions starting orthogonally from the rotation axis. 
 
In the other two types of equations, we  achieve successful answers to the existence problem of \eqref{rot}. In the case that the equation is hyperbolic at $r=0$, we obtain:
 
\begin{thm}\label{t2} If \eqref{rot} is of hyperbolic type at $r=0$, then there do not exist  radial  solutions of \eqref{rot}.  
\end{thm} 

If \eqref{rot} is parabololic, we find all solutions regardless of the intersection with the rotation axis being orthogonal or having a cusp, or even if the solution stays at a positive distance to the rotation axis.

\begin{thm}\label{t3}  If \eqref{rot} is of parabolic type, then the solutions are parametrizations of suitable circles of fixed radius.
 \end{thm}

This paper is organized as follows. In Section \eqref{sec2}, we relate the constants $a,b$ and the prescribed function $\phi$ with the character of the PDE \eqref{w1} as elliptic, hyperbolic and parabolic. We also state the character of equation \eqref{rot} at a single instant $r=r_0$. In Section \ref{sec3}, we address the existence of radial solutions of \eqref{w1} for the hyperbolic and parabolic cases. First, we prove that if the equation is hyperbolic, there are not solutions of \eqref{rot} intersecting orthogonally the rotation axis (Theorem \ref{t2}). Second, in the parabolic case we find all solutions forming all them a uniparametric family of circles of the same radius (Theorem \ref{t3}). Finally, in Section \ref{sec4} we focus on the elliptic case. We exhibit an affirmative solution to the existence problem of \eqref{rot}, proving Theorem \ref{t1}. Then, we prove uniqueness and symmetry results of the solutions of the Dirichlet problem \eqref{w1}.

\section{Types of Weingarten equations}\label{sec2}
Let us write \eqref{w1} in nonparametric form. Let $u=u(x,y)$ and suppose that $u$ satisfies  \eqref{w1}. If we define the functional 
$$
\mathfrak{F}(p,q,r,s,t)=a\frac{(1+p^2)s-2pqt+(1+q^2)r}{(1+p^2+q^2)^{3/2}}+b\frac{rs-t^2}{(1+p^2+q^2)^2}-\phi\left(\frac{1}{\sqrt{1+p^2+q^2}}\right).
$$
then    $\mathfrak{F}(u_x,u_y,u_{xx},u_{yy},u_{xy})=0$. Furthermore, the determinant of the coefficients of second order is $\mathfrak{F}_r\mathfrak{F}_s-\frac{1}{4}\mathfrak{F}_t^2=(1+p^2+q^2)^2(a^2+b\phi)$. Thus depending on the sign on the left-hand, we have an EDP of elliptic, parabolic or hyperbolic type.
  Bearing this in mind, the following definition arises: 

\begin{definition} Let $S$ be a surface satisfying \eqref{w2}.  
\begin{itemize}
\item If $a^2+b\phi>0$, the surface is    of \emph{elliptic} type.
\item If $a^2+b\phi=0$, the surface is of  \emph{parabolic} type.
\item If $a^2+b\phi<0$, the surface is of \emph{hyperbolic} type.
\end{itemize}
\end{definition}

\begin{rmk} The sphere of radius $r>0$ satisfies \eqref{w2} for different values of $a,b$ and choices of $\phi$. Indeed, the left-hand side of \eqref{w2} is $(2ar+b)/r^2$. Taking   $\phi$ the constant function $\phi=(2ar+b)/r^2$, then
$$a^2+b\phi=\frac{(ar+b)^2}{r^2}.$$
Thus the sphere satisfies \eqref{w2} for many values of $a,b$,  being elliptic or parabolic depending if   $ar+b\not=0$ or $ar+b=0$, respectively. 
\end{rmk}

We emphasize that for fixed $a,b\in\R$ and $\phi$, a given surface may have points of the three types, depending on the height of the parallel in $\S^2$ where the Gauss map $N$ lies, and eventually on the value of $\phi(\langle N,(0,0,1)\rangle)$.

Some of the results achieved in this paper only depend on the local character of the PDE \eqref{w1} as elliptic, hyperbolic or parabolic. For instance, as proved in subsequent sections, the existence of solutions of \eqref{rot} solely depends on the sign of the quantity $a^2+b\phi(1)$. For other results, as the ones exhibited in Section \ref{sec4}, the elliptic condition $a^2+b\phi>0$ must be everywhere fulfilled.

Taking into account these discussions, we settle the notation in the following definition.

\begin{definition}
Let be $u=u(r)$ a solution of \eqref{rot} and $r_0\geq0$. We say that \eqref{rot} is:
\begin{itemize}
\item of elliptic type at $r=r_0$ if $a^2+b\phi\left(\frac{1}{\sqrt{1+u'(r_0)^2}}\right)>0$;
\item of parabolic type at $r=r_0$ if $a^2+b\phi\left(\frac{1}{\sqrt{1+u'(r_0)^2}}\right)=0$; and
\item of hyperbolic type at $r=r_0$ if $a^2+b\phi\left(\frac{1}{\sqrt{1+u'(r_0)^2}}\right)<0$.
\end{itemize}
\end{definition}
The conditions of being elliptic and hyperbolic are \emph{open}, in the sense that if \eqref{rot} is elliptic or hyperbolic at some $r_0$, there exists some $\varepsilon>0$ such that \eqref{rot} is elliptic or hyperbolic for every $r\in(r_0-\varepsilon,r_0+\varepsilon)\cap[0,\infty)$.

We will simply say that \eqref{rot} is of elliptic type if $a^2+b\phi>0$ for every possible value in the argument of $\phi$, and similarly to the hyperbolic and parabolic types.

\section{Radial solutions: hyperbolic and parabolic type}\label{sec3}

In this section we investigate the existence of classical radial solutions of \eqref{rot} in the hyperbolic and parabolic cases. We first prove Theorem \ref{t2}, which is formulated again for the reader's convenience.
\begin{thm}[hyperbolic type]
If \eqref{rot} is of hyperbolic type at $r=0$, then there are not solutions of the initial value problem \eqref{rot}.
\end{thm}

\begin{pf} By contradiction, suppose that  $u=u(r)$ is a solution of the initial value problem \eqref{rot}.  
Taking limits in  \eqref{rot} as $r\rightarrow0$ and applying the L'H\^{o}pital rule to the quotient  $u'(r)/r$, we have 
$$
2au''(0)+bu''(0)^2=\phi(1),
$$
because $u'(0)=0$ in $\phi(1/\sqrt{1+u'^2})$. However, the discriminant of this equation on $u''(0)$ is  $a^2+b\phi(1)$ which is negative, obtaining a contradiction. 
\end{pf}

Now we address the existence problem of \eqref{rot} in the parabolic case. Since $a^2+b\phi=0$ everywhere, then $\phi$ is a constant function, say,  $\phi=c$. Because $a^2+bc=0$, in particular $b\not=0$. If we divide \eqref{rot} by $-b$, we can assume that $b=-1$ and $c=a^2$. Furthermore, after a change of orientation on the surface, if necessary, we can suppose that $a>0$. Note that this change of the orientation does not affect the right-hand side of \eqref{w1}, since $\phi=c$. Definitively, the class of parabolic equations \eqref{w2} reduces to the linear Weingarten relation $2aH-K=a^2$ with $a>0$. 

\begin{figure}[hbtp]
\begin{center}
\includegraphics[width=.6\textwidth]{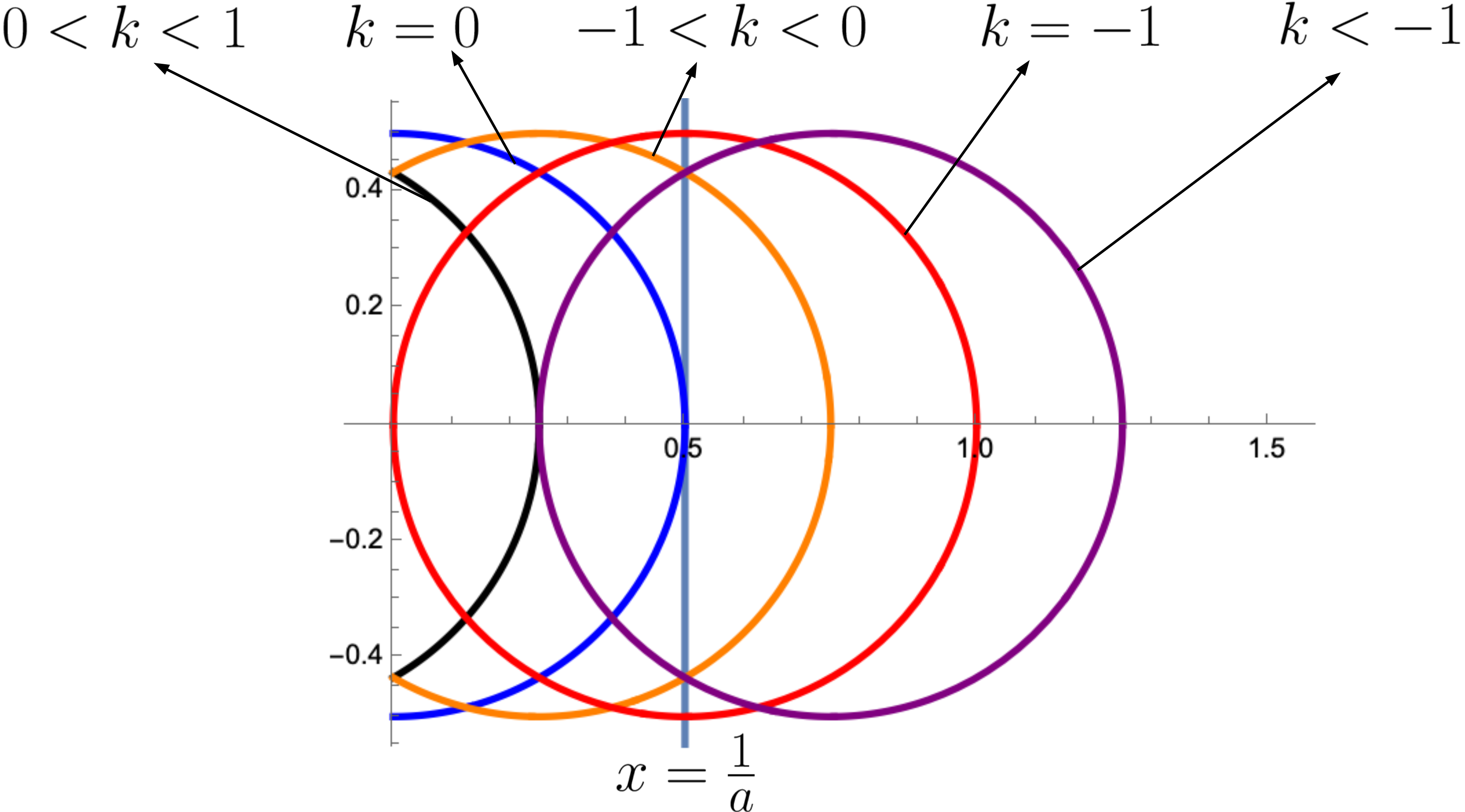} 
\end{center}
\caption{Radial solutions of the parabolic Weingarten equation.}\label{fig1}
\end{figure}

We formulate again Theorem \ref{t3} and prove it, by finding all radial solutions of \eqref{w1} independently if the surface meets or not the rotation axis. See Figure \ref{fig1}.

\begin{thm}[parabolic type] The   solutions of  
\begin{equation}\label{p0}
a\left(\frac{u''}{(1+u'^2)^{3/2}}+\frac{u'}{r\sqrt{1+u'^2}}\right)-\frac{u''u'}{r(1+u'^2)^2}=a^2, \quad a>0,
\end{equation}
are circles of radius $a$.
\end{thm}

\begin{pf}
From \eqref{p0},  
$$
\frac{u''}{(1+u'^2)^{3/2}}\left( a-\frac{u'}{x\sqrt{1+u'^2}}\right)=a\left( a-\frac{u'}{x\sqrt{1+u'^2}}\right).
$$
This implies the discussion of  two cases.
\begin{enumerate}
\item Suppose that there is $r_0>0$ such that  
$$
a\not=\frac{u'(r_0)}{x\sqrt{1+u'(r_0)^2}}.
$$
Then in an interval around $r_0$, 
$$
\frac{u''}{(1+u'^2)^{3/2}}=a.$$
Then it is immediate that  
\begin{equation}\label{p1}
u(r)=\pm\frac{1}{a}\sqrt{1-(ar+k)^2}+m, 
\end{equation}
for some constants $ k,m\in\R$. It is straightforward that $u$ parametrizes a circle of radius $1/a$.

\item Suppose  
$$
a =\frac{u'(r)}{x\sqrt{1+u'(r)^2}}
$$
for all $r$. Solving this equation,  
$$u(r)=\pm\frac{1}{a}\sqrt{1-a^2r^2}+m,\quad m\in\R. $$
Then $u$   parametrizes a circle   centered at $r=0$ of radius $1/a$. Let us notice that this solutions is particular of \eqref{p1} with $k=0$.
\end{enumerate}

\end{pf}

Studying in detail each choice of $k$ in \eqref{p1}, we conclude the next classification of the radial solutions of \eqref{w2} in the parabolic case. 
\begin{cor}\label{cor-pa}
 The radial solutions of  the equation 
$$a\, \mbox{div} \dfrac{Du}{\sqrt{1+|Du|^2}} -  \dfrac{\mbox{det} D^2u}{(1+|Du|^2)^2}=a^2$$
are:

\begin{enumerate}
\item The vertical straight-line at $r_0=1/a$. 
\item From the solutions of \eqref{p1}, the constant $k$ must be less than $1$. Furthermore, 
\begin{enumerate}
\item If $k\in(0,1)$ we obtain a one-parameter family of minor subarcs of the circle of radius $1/a$ that intersect the $z$-axis at two cusp points.
\item If $k=0$ we obtain a half-circle centered at the $z$-axis of radius $1/a$.  
\item If $k\in(-1,0)$ we obtain a one-parameter family of major subarcs of the circle of radius $1/a$ that intersect the $z$-axis at two cusp points.
\item If $k=-1$ we obtain the full circle of radius $1/a$ intersecting tangentially the $z$-axis.
\item If $k<-1$ we obtain the full circle of radius $1/a$ strictly contained in the halfplane $r>0$.  
\end{enumerate}
\end{enumerate} 
\end{cor}

\begin{pf}
A particular case to consider of radial solutions occurs when the generating curve is not a graph on the $r$-axis, that is, it is a vertical straight-line at $r=r_0$. In such a case, the surface is a circular cylinder, hence $K=0$ and $H=1/(2r_0)$, obtaining the example of item (1). It only remains to notice that from the solutions of \eqref{p1}, we deduce that $k<1$ because $|ar+k|<1$ and $r>0$. In such a case, the description of item (2) is obvious by varying $k$ from $1$ to $-\infty$.  
\end{pf}

In terms of surfaces of revolution, we conclude
\begin{cor} The rotational linear Weingarten surfaces of parabolic type are circular cylinders, spheres, embedded tori of revolution and a 1-parameter family of non-complete examples intersecting the rotation axis at cusp points and whose profile curves are arcs of a circle of fixed radius.
\end{cor}


\section{Existence of radial surfaces: elliptic case}\label{sec4}
 
In this last section we study \eqref{w1} in the elliptic case. First, we prove Theorem \ref{t1}, whose formulation is stated again.

\begin{thm}\label{t-elli} If the equation in \eqref{rot} is elliptic at $r=0$, there is $R>0$ such that the initial value problem \eqref{rot} has a  solution in $[0,R]$.
\end{thm}

\begin{pf}
Multiplying   \eqref{rot} by $r$, it is immediate that we can write \eqref{rot} as  
$$
\left(\frac{ru'}{\sqrt{1+u'^2}}\right)'+\frac{b}{2a}\left(\frac{u'^2}{1+u'^2}\right)'=\frac{r}{a}\phi\left(\frac{1}{\sqrt{1+u'^2}}\right).$$
Define the functions $f,g\colon\R \to\R$ by  
$$f(y)=\frac{y}{\sqrt{1+y^2}},
\quad g(y)=\frac{1}{a}\phi\left(\frac{1}{\sqrt{1+y^2}}\right).$$
Now we write the above equation as 
$$ rf(u')+\frac{b}{2a}f(u')^2=\int_0^r tg(u'(t))\, dt.$$
After solving $f(u')$ and eventually $u'$, we define the operator
$$
(\mathsf{T}u)(r)=\int_0^rf^{-1}\left(\frac{2a}{b}\left(-s+\sqrt{s^2+\frac{b}{a}\int_0^stg(u'(t))\, dt}\right)\right)\, ds.
$$
It can be easily proved that $u$ is a solution of the problem \eqref{rot} if $u$ is a fixed point of the operator $\mathsf{T}$. In this setting, we exhibit the existence of $R>0$ such that $\mathsf{T}$ a contraction in the space $C^1([0,R])$ endowed with the usual norm $\|u\|=\|u\|_\infty+\|u'\|_\infty$. Denote $L_{f^{-1}}$ and $L_g$ the Lipschitz constants of $f^{-1}$ and $g$ in $[-\epsilon,\epsilon]$, respectively, provided $\epsilon<1$. For all $u,v\in C^1([0,R])$, we have  
$$\|\mathsf{T}u-\mathsf{T}v\|=\|\mathsf{T}u-\mathsf{T}v\|_\infty+\|(\mathsf{T}u)'-(\mathsf{T}v)'\|_\infty,$$
We first study the term   $\|\mathsf{T}u-\mathsf{T}v\|_\infty$ because $\|(\mathsf{T}u)'-(\mathsf{T}v)'\|_\infty$ is similar. Given  two functions $u,v$ in the ball $\overline{\mathcal{B}(0,\epsilon)}\subset (C^1([0,R]),\|\cdot\|)$ and for all $r\in [0,R]$, where $R$ will be determined later, we have  
\begin{equation}\label{e1}
\begin{split}
&|(\mathsf{T}u)(r)-(\mathsf{T}v)(r)|\leq\\
&
\int_0^r\Bigg|f^{-1}\left(\frac{2a}{b}\left(-s+\sqrt{s^2+\frac{b}{a}\int_0^stg(u')\, dt}\right)\right)-\\
&f^{-1}\left(\frac{2a}{b}\left(-s+\sqrt{s^2+\frac{b}{a}\int_0^stg(v')\, dt}\right)\right)\Bigg|\, ds\\
&\leq \frac{2a}{b}L_{f^{-1}}\int_0^r\Bigg|\sqrt{s^2+\frac{b}{a}\int_0^stg(u')\, dt}-\sqrt{s^2+\frac{b}{a}\int_0^stg(v')\, dt}\Bigg|\, ds,
\end{split}
\end{equation}
where $L_{f^{-1}}$ stands for the Lipschitz constant of the function $f^{-1}$. Using the L'H\^{o}pital rule, the behavior of the function   $\int_0^stg(u')\, dt$ at $s=0$ comparing with  $s^2$ is
$$
\lim_{s\rightarrow0}\frac{\int_0^stg(u')\, dt}{s^2}=  \lim_{s\rightarrow0}\frac{sg(u'(s))}{2s}=\frac{\phi(1)}{2a}.
$$

Therefore
$$\int_0^stg(u')\, dt=c_0s^2+o(s^2),\hspace{.5cm} c_0=\frac{\phi(1)}{2a}.$$
Following the argument in \eqref{e1},  
\begin{equation*}
\begin{split}
&|(\mathsf{T}u)(r)-(\mathsf{T}v)(r)|\\
&\leq \frac{2a}{b}L_{f^{-1}}\int_0^r\Bigg|\sqrt{s^2+\frac{b}{a}\int_0^stg(u')\, dt}-\sqrt{s^2+\frac{b}{a}\int_0^stg(v')\, dt}\Bigg|\, ds\\
&=2L_{f^{-1}}\int_0^r\frac{|\int_0^st(g(u')-g(v'))\, dt|}{\sqrt{s^2+\frac{b}{a}\int_0^stg(u')\, dt}+\sqrt{s^2+\frac{b}{a}\int_0^stg(v')\, dt}}\, ds\\
&\leq 2L_{f^{-1}}\int_0^r\frac{\int_0^st|g(u')-g(v')|\, dt}{\sqrt{s^2+\frac{b}{a}\int_0^stg(u')\, dt}+\sqrt{s^2+\frac{b}{a}\int_0^stg(v')\, dt}}\, ds.
\end{split}
\end{equation*}
Now, using the Lipschitz constant   $L_g$  of $g$, we have
\begin{equation*}
\begin{split}
&\leq 2L_{f^{-1}}L_g\int_0^r\frac{\int_0^st|u'(t)-v'(t)|\, dt}{\sqrt{s^2+\frac{b}{a}\int_0^stg(u')\, dt}+\sqrt{s^2+\frac{b}{a}\int_0^stg(v')\, dt}}\, ds\\
&\leq L_{f^{-1}}L_g\|u-v\|\int_0^r\frac{s^2}{\sqrt{s^2+\frac{b}{a}\int_0^stg(u')\, dt}+\sqrt{s^2+\frac{b}{a}\int_0^stg(v')\, dt}}\, ds\\
&\leq L_{f^{-1}}L_g\|u-v\|\int_0^r\frac{s}{\sqrt{1+\frac{b}{a}\frac{\int_0^stg(u')\, dt}{s^2}}+\sqrt{1+\frac{b}{a}\frac{\int_0^stg(v')\, dt}{s^2}}}\, ds.
\end{split}
\end{equation*}
Taking into account that  $\int_0^stg(u')\, dt=c_0s^2+o(s^2)$, we follow the above expression: 
\begin{equation*}
\begin{split}
&=   L_{f^{-1}}L_g\|u-v\|\int_0^r\frac{s\, ds}{\sqrt{1+\frac{b}{a}(c_0+o(1))}+\sqrt{1+\frac{b}{a}(c_0+o(1))}}.
\end{split}
\end{equation*}
Bearing in mind that since $c_0=\phi(1)/(2a)$ and
$$
1+\frac{b}{a}c_0=\frac{1}{2}+\frac{a^2+b\phi(1)}{2a^2}>\frac{1}{2}>0,
$$
we conclude that for $r$ close to $0$, the denominator in the above integral can be upper bounded by a constant $C>0$. Hence,  
$$|(\mathsf{T}u)(r)-(\mathsf{T}v)(r)|\leq CL_{f^{-1}}L_g\|u-v\|\int_0^rs\, ds=CL_{f^{-1}}L_g\frac{r^2}{2}\|u-v\|.$$
Let $R_1$ be sufficiently small such that the constant  $K_1=CL_{f^{-1}}L_gR_1^2/2$ is less than $1$. Then $\|\mathsf{T}u-\mathsf{T}v\|_{\infty}\leq K_1\|u-v\|$.  

As we have said, a similar argument works with     $\|(\mathsf{T}u)'-(\mathsf{T}v)'\|_\infty$. In virtue of the definition of $\mathsf{T}$,
$$
(\mathsf{T}u)'(r)=f^{-1}\left(\frac{2a}{b}\left(-r+\sqrt{r^2+\frac{b}{a}\int_0^rtg(u'(t))\, dt}\right)\right).$$
Thus
\begin{equation*}
\begin{split}
&|(\mathsf{T}u)'(r)-(\mathsf{T}v)'(r)|\\
&\leq L_{f^{-1}}\frac{2a}{b}\left|\sqrt{r^2+\frac{b}{a}\int_0^rtg(u'(t))\, dt}-\sqrt{r^2+\frac{b}{a}\int_0^rtg(v'(t))\, dt}\right|\\
&\leq  L_{f^{-1}}\frac{2\int_0^rt|g(u'(t))-g(v'(t))|dt}{\sqrt{r^2+\frac{b}{a}\int_0^rtg(u'(t))\, dt}+\sqrt{r^2+\frac{b}{a}\int_0^rtg(v'(t))\, dt}}.
\end{split}
\end{equation*}
Again, since $g$ is Lipschitz,
\begin{equation*}
\begin{split}
&\leq L_{f^{-1}}L_g\|u-v\|\frac{r^2}{\sqrt{r^2+\frac{b}{a}\int_0^rtg(u'(t))dt}+\sqrt{r^2+\frac{b}{a}\int_0^rtg(v'(t))dt}}\\
&=L_{f^{-1}}L_g\|u-v\|\frac{r}{\sqrt{1+\frac{b}{a}\frac{\int_0^rtg(u'(t))\, dt}{r^2}}+\sqrt{1+\frac{b}{a}\frac{\int_0^rtg(v'(t))\, dt}{r^2}}}\\
&\leq L_{f^{-1}}L_g\|u-v\|\frac{r}{\sqrt{1+\frac{b}{a}(c_0+o(1))}+\sqrt{1+\frac{b}{a}(c_0+o(1))}}\\
&\leq  CL_{f^{-1}}L_g\|u-v\|r.
\end{split}
\end{equation*}
Let $R_2>0$ be sufficiently small so the constant $K_2=CL_{f^{-1}}L_gR_2$ is less than $1$. With this constant, if $r\in [0,R_2]$, we have 
$\|(\mathsf{T}u)'-(\mathsf{T}v)'\|\leq K_2\|u-v\|$. Finally, we now choose the value $R$ as  $R=\min\{R_1,R_2\}$. Thus if $r\in [0,R]$, we have
$$\|\mathsf{T}u-\mathsf{T}v\|<\max\{K_1,K_2\}\|u- v\|,$$
proving that   the operator $\mathsf{T}$ is contractible in   $C^1([0,R])$. This proves the existence of  a fixed point $u\in C^1([0,R])\cap C^2((0,R])$.

Finally, we prove that the solution $u$ extends with $C^2$-regularity at $r=0$. By taking limits as $r\to 0$,  and by L'H\^{o}pital rule again on the quotient $u'(r)/r$, we conclude  
$$2au''(0)+bu''(0)^2=\phi(1).$$
In this case, 
$$u''(0)=\frac{-a\pm\sqrt{a^2+b\phi(1)}}{b},$$
which has a solution by the elliptic condition  $a^2+b\phi(1)>0$.
\end{pf}

\begin{rmk}
In the proof of Theorem \ref{t-elli} we can relax the $C^1$-regularity of $\phi(y)$ to just being Lipschitz continuous around $y=1$. Indeed, all the arguments of the proof are local and one of the hypotheses needed is $g(y)=\frac{1}{a}\phi\left(\frac{1}{\sqrt{1+y^2}}\right)$ to be Lipschitz around $y=0$, i.e. $\phi(y)$ to be Lipschitz around $y=1$.
\end{rmk}

In the remaining of this paper we assume that the equation in \eqref{w1} is elliptic, i.e. $a^2+b\phi>0$ for every possible argument of the function $\phi$. This global elliptic condition will allow us to obtain results concerning the global behavior of the solutions of equation \eqref{w1}. 

First, we prove the uniqueness of the Dirichlet problem \eqref{w1}. Here we use the  comparison principle for fully nonlinear elliptic PDEs (\cite[Th. 17.1]{gt}) to the functional $\mathfrak{F}$ which, in our context of Weingarten surfaces, asserts that if $u_1$ and $u_2$ are two functions defined in $\Omega$ such that $\mathfrak{F}[u_1]\geq\mathfrak{F}[u_2]$ in $\Omega$ and $u_1\leq u_2$ on $\partial\Omega$, then $u_1\leq u_2$ in $\Omega$. Furthermore, if $\mathfrak{F}[u_1]>\mathfrak{F}[u_2]$ in $\Omega$, then $u_1<u_2$ in $\Omega$. Similarly the functional $\mathfrak{F}$ satisfies a maximum principle in the sense that if $\mathfrak{F}[u_1]=\mathfrak{F}[u_2]$, $u_1=u_2$ at some point $x_0\in\Omega$ and $u_1\geq u_2$ in an open set of $x_0$, then $u_1=u_2$ in that open set. 

We now prove the uniqueness of the Dirichlet problem \eqref{w1} assuming arbitrary continuous boundary values. 
 
\begin{proposition}\label{pr-uni}
Suppose that the equation in \eqref{w1} is elliptic. If the Dirichlet problem \eqref{w1} has a solution for continuous boundary values $u=\varphi$ on $\partial\Omega$, then the solution  is unique.
\end{proposition}

\begin{pf} The argument is standard using the comparison principle and the fact  the vertical translations of $\R^3$ are isometries that preserve the solutions of \eqref{w1}. If $u_1$ and $u_2$ are two such solutions, we move the graph $S_1$ of  $u_1$ downwards until that it does not intersect $S_2$, the graph of $u_2$. This is possible because $S_1$ and $S_2$ are compact surfaces. Now we move $S_1$ upwards until reaching a first contact with $S_2$. Then the (interior or boundary version of the) maximum principle asserts that $S_1=S_2$, that is, $u_1=u_2$ in $\Omega$.
\end{pf}

We finish this section proving that in case that the equation \eqref{w1} is elliptic, the solutions of the Dirichlet problem \eqref{w1} are radial if $\Omega$ is an Euclidean ball. The method comes back to the known technique of moving planes. First, we need the next result, which has its own interest.

\begin{proposition}\label{pr-no}
Assume that the equation \eqref{w1} is elliptic. If the Dirichlet problem \eqref{w1} has a solution $u$,  then $u$ has constant sign in $\Omega$.
\end{proposition}

\begin{pf} By contradiction, suppose that $u$ changes of sign in $\Omega$.   Let  $x_0,x_1\in\Omega$ be the points where $u$ attain its minimum and maximum and suppose without loss of generality that  $u(x_0)\leq 0<u(x_1)$. In particular, $Du(x_0)=Du(x_1)=0$. Let $v_0,v_1\colon\Omega\to\R$ be the constant functions defined by $v_0(x,y)=u(x_0)$ and $v_1(x,y)=u(x_1)$. Since $v_0\leq u$ in a neighborhood of $x_0$ and $u\leq v_1$ around $x_1$, and because  the functional $\mathfrak{F}$ is elliptic, the comparison principle implies
$$\mathfrak{F}[v_0]< \mathfrak{F}[u],\quad \mathfrak{F}[u]< \mathfrak{F}[v_1].$$
Since $\mathfrak{F}[u]=0$ and $\mathfrak{F}[v_0]=\mathfrak{F}[v_1]=-\phi(1)$, we obtain  a contradiction.
\end{pf}

Once proved that the solution has sign in $\Omega$, or equivalently, the surface that determines lies completely at one side of the coordinate plane $z=0$, we can prove that if $\Omega$ is a round disc, then the solution of \eqref{w1} is radially symmetric, or equivalently, the surface is rotational about the $z$-axis.  Here we follow the moving plane method of Alexandrov (\cite{al}), see also \cite{gi,se}. The arguments are standard and the key issue is that the equation is elliptic and Proposition \ref{pr-no}. We give a brief proof, stating the result in its more general assumption of the domain $\Omega$. In the next result, we will denote by $(x_1,x_2,x_3)$ the coordinates of $\R^3$.

\begin{thm} Suppose that $\Omega\subset\R^2$ is a bounded smooth domain, convex in the $x$-direction and symmetric about the line $x_1=0$. If \eqref{w1} is of elliptic type, then any solution $u\in C^2(\overline{\Omega})$ of \eqref{w1} with Dirichlet condition $u=0$ along $\partial\Omega$ is also symmetric about the line $x_1=0$. 
\end{thm}

\begin{pf} From Proposition \ref{pr-no} we know that $u$ has constant sign. Without loss of generality, we suppose that $u<0$ in $\Omega$. Since the function $\phi$ depends on $1/\sqrt{1+|Du|^2}$ (or $\phi=\phi(\langle N,v\rangle)$ in \eqref{w2}), then the reflections about a vertical plane of a surface that satisfies the Weingarten equation \eqref{w2} are surfaces satisfying the same equation. 

For $t\leq 0$, let $\Omega_t=\Omega\cap \{x_1\leq t\}$. If $A\subset \R^2$, with the notation $A^*$ we stand for the reflection of $A$  about the line of equation $x_1=t$, that is,   $A^*=\{(2t-x_1,x_2):(x_1,x_2)\in A\}$. Define on $\Omega_t^*$ the function $u_t$ obtained by reflection about the line $x_1=t$, $u_t(x_1,x_2)=u(2t-x_1,x_2)$. Then $u_t$ satisfies \eqref{w1} in $\Omega_t^*$. We begin with the method of moving planes doing reflection about the line $x_1=t$ for $t$ near $-\infty$. Since $\Omega$ is bounded and after the first touching point $t_1<0$ with $\partial\Omega$, we have $u<u_t$ in $\Omega_t$ for $t\in (t_1,t_1+\epsilon)$ for some $\epsilon>0$ sufficiently small. Moving $t\nearrow 0$, and by the compactness of $\Omega$, let 
$$t_0=\sup\{t<0:u<u_t \mbox{ in }\Omega_t^*\}.$$
If $t_0<0$ and because $(\partial\Omega_t\cap\partial\Omega)^*\subset\Omega$, $u<0$ in $\Omega$ and the convexity of $\Omega$ in the $x_1$-direction, there is $x_0\in\Omega_t^*$ such that $u=u_{t_0}$ at $x_0$. Since $u\leq u_{t_0}$ and $\mathfrak{F}[u]=\mathfrak{F}[u_{t_0}]$ in $\Omega_{t_0}^{*}$, then $u=u_{t_0}$ by the maximum principle. Using an argument of connectedness, this implies that the line $x_1=t_0$ is a line of symmetry of $u$, which it is false because $\Omega_{t_0}^*\cup\Omega_{t_0}\not=\Omega$. 

Thus $t_0=0$ and $u<u_t$ in $\Omega_t^*$ for all $t<0$. By the symmetry of $\Omega$ with respect to the line $x_1=0$, there is $x_0\in \partial\Omega_0^*\cap\partial \Omega$ such that $u(x_0)=u_0(x_0)=0$. Using the maximum principle of elliptic equations in its boundary version, we conclude that that $u=u_0$ in $\Omega\cap\Omega_0^*$, proving the result.
\end{pf}

As  consequence of this theorem, together with Theorem \ref{t-elli} and Proposition \ref{pr-uni}, we conclude the following consequence.
\begin{cor} 
Assume that the equation \eqref{w1} is elliptic. Then there is $R>0$ such that the Dirichlet problem \eqref{w1} in the ball $B(0,R)$ has a unique solution. Moreover, this solution is radial.
\end{cor}

\section*{Acknowledgments} Rafael L\'opez belongs to the Project  I+D+i PID2020-117868GB-I00, supported by MCIN/ AEI/10.13039/501100011033/


\begin{thebibliography}{00}

\bibitem{al}  A. D. Alexandrov,  A characteristic property of spheres. Ann. Mat. Pura Appl. 58 (1962),  303--315.
\bibitem{an} B. Andrews, Contraction of convex hypersurfaces in Euclidean space. Calc. Var. Partial Differ. Eq. 2  (1994), 151--171.
\bibitem{bss} A. Barros, J. Silva, P. Sousa, Rotational linear Weingarten surfaces into the Euclidean sphere. Israel J. Math. 192 (2012), 819--830. 

\bibitem{bgm} A. Bueno, J. A. G\'alvez, P. Mira, The global geometry of surfaces with prescribed mean curvature in $R^3$. Trans. Amer. Math. Soc. 373 (2020),  4437--4467.

\bibitem{bgm2} A. Bueno, J. A. G\'alvez, P. Mira, Rotational hypersurfaces of prescribed mean curvature. J. Differ. Eq. 268  (2020), 2394--2413.
\bibitem{gm} J. A. G\'alvez, P.  Mira, Rotational symmetry of Weingarten spheres in homogeneous three-manifolds.  J. Reine Angew. Math. 773 (2021), 21--66. 

\bibitem{gi} B. Gidas,W. M.  Ni, L. Nirenberg,   Symmetry and related properties via the maximum principle. Commun. Math. Phys. 68  (1979), 209--243.

\bibitem{gt} D. Gilbarg, N. Trudinger, Elliptic Partial Differential Equations of Second Order, (Classics in Mathematics). Springer, 2001.

\bibitem{hsi}   G. Huisken,  , C.  Sinestrari,     Convexity estimates for mean curvature flow   and singularities of mean convex surfaces.  Acta Math.  183  (1999), 45--70.


\bibitem{il}  T. Ilmanen,   Elliptic regularization and partial regularity for motion by mean curvature. Mem. Amer. Math. Soc. 108, x+90 (1994).
  
  
 \bibitem{lo1} R. L\'opez, Rotational linear Weingarten surfaces of hyperbolic type. Israel J. Math. 167 (2008), 283--301. 
 
\bibitem{lo2} R. L\'opez,  Parabolic Weingarten surfaces in hyperbolic space. Publ. Math. Debrecen 74  (2009), 59--80.
 
\bibitem{mr} F. Morabito, M. M. Rodr\'{\i}guez, Classification of rotational special Weingarten surfaces of minimal type in $H^2\times R$ and $S^2\times R$. Math. Z. 273 (2013),   379--399.

\bibitem{oo} F. Obersnel, P. Omari, Revisiting the sub-and super-solution method for the classical radial solutions of the mean curvature equation. Open Math. 18 (2020),  1185--1205.

\bibitem{rs} H. Rosenberg, R. Sa Earp, The geometry of properly embedded special surfaces in $\R^3$, e.g., surfaces satisfying $aH+bK=1$, where $a$ and $b$ are positive. Duke Math. J., 73 (1994), 291--306.  
 
\bibitem{st} R. Sa Earp et E. Toubiana, Classification des surfaces de type Delaunay et applications. Amer. J. Math. 221 (1999), 671--700.
\bibitem{se} J. Serrin,   A symmetry problem in potential theory. Arch. Ration. Mech. Anal. 43 (1971), 304--318

\bibitem{ur} J. Urbas,   An expansion of convex hypersurfaces. J. Differential Geom. 33 (1991), 91--125.
\end{thebibliography}
\end{document}